\documentclass[12pt]{amsart}
\usepackage{upref}
\usepackage[latin1]{inputenc}
\usepackage{mathptmx}

\newcommand{\aaf}{\mathfrak a}
\newcommand{\pp}{\mathfrak p}

\newcommand{\qq}{\mathfrak q}

\newcommand{\Cl}{{\rm Cl}}

\newcommand{\ini}{{\rm in}}
\newcommand{\rank}{{\rm rank}}

\newcommand{\Ker}{{\rm Ker}}

\newcommand{\GL}{{\rm GL}}
\newcommand{\SL}{{\rm SL}}

\newcommand{\Dirsum}{{\bigoplus}}

\newcommand{\pnt}{{\raise 0.5mm \hbox{\large\bf.}}}

\newcommand{\sep}{{\,|\,}}
\def\SSS{{\mathcal{S}}}

\newcommand{\sym}{{\textup{sym}}}
\newcommand{\tr}{{\textup{tr}}}
\let\epsilon=\varepsilon
\let\phi=\varphi
\let\theta=\vartheta

\newcommand{\ZZ}{\mathbb{Z}}
\newcommand{\RR}{\mathbb{R}}

\newtheorem{Theorem}{\bf Theorem}[section]
\newtheorem{Lemma}[Theorem]{\bf Lemma}
\newtheorem{Corollary}[Theorem]{\bf Corollary}
\newtheorem{Proposition}[Theorem]{\bf Proposition}

\theoremstyle{definition}
\newtheorem{Remark}[Theorem]{\bf Remark}

\newtheorem{Example and Definition}[Theorem]{\bf Example and Definition}

\title{Initial algebras of determinantal rings, Cohen-Macaulay and Ulrich ideals}

\author{Winfried Bruns}
\address{FB Mathematik/Informatik, Universit\"at Osnabr\"uck, 49069 Osnabr\"uck, Germany}
\email{winfried@mathematik.uni-osnabrueck.de}

\author{Tim R\"omer}
\address{FB Mathematik/Informatik, Universit\"at Osnabr\"uck, 49069 Osnabr\"uck, Germany}
\email{troemer@mathematik.uni-osnabrueck.de}

\author{Attila Wiebe}
\address{FB Mathematik, Universit\"at Duisburg-Essen, Campus Essen,  45117 Essen, Germany}
\email{attila.wiebe@ruhrgas.com}

\begin{document}

\begin{abstract}
We study initial algebras of determinantal rings, defined by
minors of ge\-neric matrices, with respect to their classical
generic point. This approach leads to very short proofs for the
structural properties of determinantal rings. Moreover, it allows
us to classify their Cohen-Macaulay and Ulrich ideals.
\end{abstract}
\maketitle
%
%
%
\section{Introduction}
Let $K$ be a field and $X$ an $m \times n$ matrix of indeterminates over $K$.
Let $K[X]$ denote the polynomial ring generated by all the indeterminates
$X_{ij}$. For a given positive integer $r \leq \min\{m,n\}$ we consider
the determinantal ideal $I_{r+1}=I_{r+1}(X)$ generated by all $r+1$ minors of $X$
if $r<\min\{m,n\}$ and $I_{r+1}=(0)$ otherwise.
Let $R_{r+1}=R_{r+1}(X)$ be the determinantal ring $K[X]/I_{r+1}$.

Determinantal ideals and rings are well-known objects and the study of these
objects has many connections with algebraic geometry, invariant theory, representation
theory and combinatorics. See Bruns and Vetter \cite{BV} for a detailed discussion.

In the first part of this paper we develop an approach to
determinantal rings via initial algebras. We cannot prove new
structural results on the rings $R_{r+1}$ in this way, but the
combinatorial arguments involved are extremely simple. They yield
quickly that $R_{r+1}$, with respect to its classical generic
point, has a normal semigroup algebra as its initial algebra.
Using general results about toric deformations and the properties
of normal semigroup rings, one obtains immediately that $R_{r+1}$
is normal, Cohen-Macaulay, with rational singularities in
characteristic $0$, and $F$-rational in characteristic $p$.

Toric deformations of determinantal rings have been constructed by
Sturmfels \cite{Stu} for the coordinate rings of Grassmannians
(via initial algebras) and Gonciulea and Lakshmibai \cite{GL} for
the class of rings considered by us. The advantage of our
approach, compared to that of \cite{GL}, is its simplicity.

Moreover, it allows us to determine the Cohen-Macaulay and Ulrich
ideals of $R_{r+1}$. Suppose that $1\leq r<\min\{m,n\}$ and let
$\pp$ (resp.\ $\qq$) be the ideal of $R_{r+1}$ generated by the
$r$-minors of the first $r$ rows (resp.\ the first $r$ columns) of
the matrix $X$. The ideals $\pp$ and $\qq$ are prime ideals of
height one and hence they are divisorial, because $R_{r+1}$ is a
normal domain. The divisor class group $\Cl(R_{r+1})$ is
isomorphic to $\ZZ$ and is generated by the class $[\pp]=-[\qq]$.
(See Bruns and Herzog \cite[Section 7.3]{BH} or
\cite[Section 8]{BV}.) The symbolic powers of $\pp$ and $\qq$ coincide with the
ordinary ones. Therefore the ideals $\pp^k$ and $\qq^k$ represent
all reflexive rank $1$ modules. The goal of the last section is to
show that $\pp^k$ (resp.\ $\qq^k$) is a Cohen-Macaulay ideal if
and only if $k\leq m-r$ (resp.\ $k\leq n-r$). In addition we prove
that the powers $\pp^{m-r}$ and $\qq^{n-r}$ are even Ulrich
ideals.

%
%
%
\section{Standard bitableaux}

Let $K$ be a field. For the study of the determinantal rings
$R_{r+1}$ we use the approach of standard bitableaux for which one
considers all minors of the matrix $X$ as generators for the
$K$-algebra $K[X]$ and not only the 1-minors $X_{ij}$. Hence
products of minors appear as ``monomials''.

Let $1\leq t \leq \min\{m,n\}$. Denote the determinant of the
matrix $X'=(X_{a_ib_j}\colon i=1,\ldots,t,\allowbreak\
j=1,\ldots,t )$ by
$$
[a_1\ldots a_t|b_1\ldots b_t].
$$
We require that $1\leq a_1<\cdots<a_t \leq m$ and $1\leq
b_1<\cdots<b_t \leq n$. We call $[a_1\ldots a_t|b_1\ldots b_t]$ a
\emph{minor} of $X$ and $t$ its \emph{size}. A \emph{bitableau}
$\Delta$ is a product of minors
$$
\prod_{i=1}^w [a_{i1}\ldots a_{it_i}|b_{i1}\ldots b_{it_i}]\quad
\text{ such that } t_1\geq \cdots \geq t_w.
$$
By convention the value of the empty minor $[\,|\,]$ is 1. The
\emph{shape} of $\Delta$ is the sequence $(t_1,\ldots,t_w)$. The
name bitableau is motivated by the graphical description of
$\Delta$ as a pair of so-called Young tableaux, and we will also
write $\Delta=(a_{ij}\sep b_{ij})$. We consider a partial order on
the set of all bitableau:
\begin{multline*}
[a_{1}\ldots a_{t}|b_{1}\ldots b_{t}] \preceq [c_{1}\ldots
c_{u}|d_{1}\ldots d_{u}]\\
\iff\quad t \geq u \text{ and } a_i \leq c_i,\ b_i \leq d_i,\
i=1,\ldots,u.
\end{multline*}
A product $\Delta=\delta_1\cdots\delta_w$ of minors
$\delta_i=[a_{i1}\ldots a_{it_i}|b_{i1}\ldots b_{it_i}]$
is a \it standard bitableau \rm if
$$
\delta_1 \preceq \cdots \preceq \delta_w,
$$
i.\ e. in each ``column'' of the bitableau the indices are
non-decreasing from top to the bottom. 
(The empty product is also standard.)
The letter $\Sigma$ is reserved for standard bitableaux. 
The fundamental straightening law of
Doubilet-Rota-Stein \cite{DRS} says that every element of $K[X]$
has a unique presentation as a $K$-linear combination of standard
bitableaux. Hence these elements are a $K$-vector space basis of $K[X]$
and $K[X]$ is an algebra with straightening law (ASL for
short) on the set of standard bitableaux. See \cite{BV} or Bruns
and Conca \cite{BC4} for a detailed introduction.

We let $\SSS_r$ denote the set of all standard bitableaux whose left
tableau has entries in $\{1,\dots,m\}$, whose right tableau has
entries in $\{1,\dots,n\}$, and whose shape $(s_1,\dots,s_u)$ is
bounded by the condition $s_1\le r$.

For a (standard) bitableau $\Sigma$ and an $m\times n$ matrix
$A=(a_{ij})$ over some $K$-algebra $B$ we let $\Sigma_A$ denote the image
of $\Sigma$ under the homomorphism $K[X]\to B$ defined by the
substitution $X_{ij}\mapsto a_{ij}$. However, for simplicity we
will not explicitly indicate the passage from $K[X]$ to its
residue class ring $R_{r+1}$.

\begin{Theorem}\label{Hodge}
The (residue classes of the) standard bitableaux $\Sigma\in\SSS_r$
generate $R_{r+1}$ as a vector space over $K$.
\end{Theorem}

The proof of this theorem, essentially due to Hodge, is to be found
in many sources. It is most easily proved by dehomogenization of
its companion result for the subalgebra of $K[X]$ spanned by the
maximal minors; for example, see \cite{BV}.
%
%
%
\section{Initial algebras}
The classical ``generic point'' for $R_{r+1}$ is the homomorphism
$$
\phi:R_{r+1}\to K[Y,Z]
$$
where $Y$ is an $m\times r$ matrix of indeterminates, $Z$ is an
$r\times n$ matrix of indeterminates, and the homomorphism is induced by the
substitution of the $(i,j)$-th entry of the product $YZ$ for
$X_{ij}$. The homomorphism $K[X]\to K[Y,Z]$ factors through
$R_{r+1}$ since $\rank (YZ)=r$.

On $K[Y,Z]$ we introduce a term order by first listing the
variables of $Y$ column by column from bottom to top, starting
with the first column, and then the entries of $Z$ row by row from
right to left:
$$
Y_{m1}>Y_{m-11}>\dots >Y_{11}>Y_{m2}>\dots >Y_{1r}>
Z_{1n}> \dots > Z_{11} > Z_{2n}>\dots >Z_{r1}.
$$
This total order is then extended to the induced degree reverse
lexicographic order on $K[Y,Z]$. Note that the restrictions of the
term orders to $K[Y]$ and $K[Z]$ are diagonal: the initial term of
a minor of $Y$ or $Z$ is the product of its main diagonal
elements. But also the initial monomials of the minors of $YZ$ are
easily found:

\begin{Lemma}\label{ini}
Let $1\leq t \leq r$.
The initial monomial of the minor $[a_1\dots a_t\sep b_1\dots
b_t]_{YZ}$ is the monomial $Y_{a_11}\cdots Y_{a_t t}Z_{1b_1}\cdots
Z_{tb_t}$.
\end{Lemma}

\begin{proof}
Suppose first that $t=r$. Then the matrix $X'=(X_{a_i b_j})$ is
the product of $Y'=(Y_{a_ij})$ and $Z'=(Z_{ib_j})$. Clearly
$$
\ini(\det(X'))=\ini(\det(Y'Z'))= \ini(\det(Y')\det(Z'))=
\ini(\det(Y'))\ini(\det(Z')),
$$
and the last term is the product of the main diagonals, as pointed
out above.

Let now $t<r$. Since we have chosen the reverse lexicographic term
order, we may delete all monomials from $[a_1\dots a_t\sep
b_1\dots b_t]_{YZ}$ that involve an indeterminate $Z_{ij}$ with
$i>t$ without loosing the initial monomial, provided at least one
term survives. But this is clearly the case: under the
substitution $Z_{ij}\mapsto 0$ for $i>t$
the minor $[a_1\dots a_t\sep
b_1\dots b_t]_{YZ}$ goes to the minor $[a_1\dots a_t\sep b_1\dots
b_t]_{\overline Y\overline Z}$ where $\overline Y$ consists of the
first $t$ columns of $Y$ and $\overline Z$ consists of the first
$t$ rows of $Z$. Now we have reached the case of maximal minors
discussed above.
\end{proof}

\begin{Proposition}\label{iniSt}\leavevmode\nopagebreak
\begin{itemize}
\item[(a)] The initial monomial of the standard bitableau
$\Sigma=(a_{ij}\sep b_{ij})$, $i=1,\dots,u$, $j=1,\dots, t_i$,
$t_1\geq\ldots\geq t_u$
is
the monomial $\prod_{i=1}^u\prod_{j=1}^{t_i}
Y_{a_{ij}j}Z_{jb_{ij}}$.

\item[(b)] If $\Sigma,\Sigma'\in\SSS_r$, $\Sigma\neq\Sigma'$, then
$\ini(\Sigma_{YZ})\neq\ini(\Sigma'_{YZ})$. In particular, the
polynomials $\Sigma_{YZ}$ are $K$-linearly independent.
\end{itemize}
\end{Proposition}

\begin{proof}
Part (a) is an immediate consequence of Lemma \ref{ini}. For part
(b) one observes that the factors $Y_{vw}$ that appear in
$\ini(\Sigma_{YZ})$ uniquely determine the $w$-th column of the
left tableau of $\Sigma$ since they indicate which indices $v$
appear in this column and determine their multiplicities. The
indices in a column are non-decreasing (from top to bottom), and
therefore the column is uniquely given by the indices and their
multiplicities. It follows that the left tableau is uniquely
determined, and a similar argument applies to the right tableau.
The linear independence follows immediately.
\end{proof}

We draw a well-known consequence.

\begin{Corollary}\label{old}
Let $K[YZ]$ denote the $K$-algebra generated by the entries of the
product matrix $YZ$.
\begin{itemize}
\item[(a)] The homomorphism $\varphi:R_{r+1}\to K[YZ]$ is an
isomorphism.
\item[(b)] The standard bitableaux $\Sigma\in\SSS_r$ form a
$K$-basis of $R_{r+1}$.
\end{itemize}
\end{Corollary}

In fact, the homomorphism maps the elements of a system of
generators of the vector space $R_{r+1}$ to a linearly independent
system in its image $K[YZ]$. In the following we will identify
$R_{r+1}$ with $K[YZ]$.

\begin{Remark}\label{eff}
(a) The above proof of the \emph{straightening law} contained in
\ref{old}(b) can be used for an effective implementation as
follows. Given an element $f\in R_{r+1}$ (so $f\in K[X]$ if
$r=\min(m,n)$), we map it to $K[YZ]$. Then the initial term of
$\phi(f)$ is determined. It determines a unique standard monomial
$\Sigma$. Next $\Sigma$ is evaluated in $R_{r+1}$ (of course, not
in $K[YZ]$!), and we replace $f$ by $f-\lambda\Sigma$ where
$\lambda$ is the leading coefficient of $\phi(f)$. Since
$f-\lambda\Sigma=0$ or $\ini(\phi(f-\lambda\Sigma))
<\ini(\phi(f))$, an iteration of the procedure must terminate
after finitely many steps.

(b) In order to avoid Theorem \ref{Hodge} in the proof of the
straightening law one would have to show that the initial monomial
of an arbitrary element in $K[YZ]$ is one of the monomials
$\ini(\Sigma_{YZ})$, $\Sigma\in\SSS_r$.

(c) If one is willing to invest the Knuth-Robinson-Schensted
correspondence, then Theorem \ref{Hodge} becomes a consequence of
Proposition \ref{iniSt}: the correspondence implies that in each
degree there exist as many standard bitableaux in $\SSS_t$,
$t=\min(m,n)$, as ordinary monomials. Together with the linear
independence of $\SSS_t$ (in whose proof Theorem \ref{Hodge} has not
been used), this implies that $\SSS_t$ is a $K$-basis of $K[X]$.
This shows \ref{Hodge} for $r=t$. The general case follows rapidly
since we have the inclusions
$$
V_{r+1}\subset I_{r+1}(X)\subset \Ker(\phi)
$$
where $V_{r+1}$ is the vector space spanned by all
$\Sigma\notin\SSS_r$ and $\SSS_r$ is mapped to a linear independent
subset of $K[YZ]$. (Note that every minor of size $>r$ is
contained in $I_{r+1}$.)

(d) We will show that the initial algebra of $R_{r+1}$ is a normal
semigroup ring. This is a direct generalization of the fact that
for $r=1$ the algebra $R_2=K[YZ]=D_2$ is a normal semigroup ring
itself.
\end{Remark}

We are in the extremely rare situation that taking initial forms
on a vector space basis is injective, and so we can immediately
describe the initial space:

\begin{Theorem}\label{main}\leavevmode\nopagebreak
\begin{itemize}
\item[(a)] The initial algebra $D_{r+1}=\ini(R_{r+1})\subset K[Y,Z]$ is
generated by the monomials $Y_{a_11}\cdots Y_{a_tt} Z_{1b_1}\cdots
Z_{tb_t}$ with $1\leq t \leq r$, $a_1<\dots< a_t$ and $b_1<\cdots<b_t$.
\item[(b)] $D_{r+1}$ is a normal semigroup ring.
\item[(c)] $R_{r+1}$ is a normal domain, Cohen-Macaulay, with
rational singularities in characteristic $0$, and $F$-rational in
characteristic $p>0$.
\end{itemize}
\end{Theorem}

\begin{proof}
(a) This is just a reformulation of Proposition \ref{iniSt}. In
fact, the subalgebra generated by the monomials given in (a) is a
$K$-vector subspace of $D_{r+1}$. On the other hand, it has the
same Hilbert function as $R_{r+1}$ (or $D_{r+1}$). This
forces equality.

(b) It is enough to show that $M^k\in D_{r+1}$ for a monomial
$M\in K[Y,Z]$ and an integer $k>0$ implies $M\in D_{r+1}$. There
exists a standard bitableau $\Sigma=(a_{ij}\sep b_{ij})$ with
$M^k=\ini(\Sigma)$. We then write $M^k$ in the form
$\prod_{i=1}^u\prod_{j=1}^{t_i} Y_{a_{ij}j}Z_{jb_{ij}}$. Since
$M^k$ is a $k$-th power and $\Sigma$ is a standard bitableau, the
first factor $\prod_{j=1}^{t_1} Y_{a_{1j}j}Z_{jb_{1j}}$ must occur
(at least) $k$ times. We split it off $M$, and conclude by
induction.

(c) follows from general theorems on flat deformation. For proofs
see \cite{BC4} or Conca, Herzog and Valla \cite{CHV}.
\end{proof}

The Cohen-Macaulay property of $R_{r+1}$ was first proved by
Hochster and Eagon \cite{HE} and the Cohen-Macaulay property of
normal semigroup rings by Hochster \cite{Ho}.

\begin{Remark}
\label{equation} For an application below we describe the set $E$
of vectors $[(\alpha_{ij}),(\beta_{uv})]\in
(\RR^{mr})\oplus(\RR^{rn})$ that appear as exponent vectors of
elements in $D_{r+1}=\ini(R_{r+1})$. It is not hard to check that
$E$ is the set of lattice points in the cone defined by the
following linear equations and inequalities:
\begin{align}
\alpha_{ij}=\beta_{uv}&=0,& j&>i,\ u>v,\label{diag1}\\
\sum_{i=j-1}^{k-1}\alpha_{ij-1}-\sum_{i=j}^k \alpha_{ij}&\ge0,&
j&=2,\dots,r,\ \ k=j,\dots,m,\\
\sum_{t=u-1}^{w-1}\beta_{u-1t}-\sum_{t=u}^w\beta_{ut}&\ge0,&
t&=2,\dots,r,\ \ w=u,\dots,n,\\
\alpha_{ij},\beta_{uv}&\ge 0,& i&>j,\ v>u\ \text{ and }\ i=j=u=v=r,\label{nonneg}\\
\sum_{i=1}^n \alpha_{ij}-\sum_{v=1}^n\beta_{jv}&=0&
j&=1,\dots,r.\label{koppel}
\end{align}
\end{Remark}

Note that for $r=\min(m,n)$ we consider an embedding of $K[X]$
into $K[Y,Z]$ which identifies the indeterminate $X_{ij}$ with the
corresponding entry of the product matrix $YZ$. Thus we can
investigate the initial ideal $\ini(J)\subset D=\ini(K[X])$ for
every ideal $J$ of $K[X]$. In particular, it is useful to consider
the ideals $I(X;\delta)$ and the residue class rings
$R(X;\delta)=K[X]/I(X;\delta)$ where $I(X;\delta)$ is generated by
all minors $\gamma\not\ge\delta$. Observe that
$R_{r+1}=R(X;\delta)$ for $\delta=[1\dots r \sep 1\dots r ]$. The
proof of the next corollary shows that we recover $D_{r+1}$ as a
retract of $D$ if we take $\delta=[1\dots r\sep 1\dots r]$.

\begin{Corollary}
Let $D$ be the initial algebra of $K[X]$. The initial ideal
$\ini(I(X;\delta))$ is a (monomial) prime ideal in $D$. Therefore
$R(X;\delta)$ is a normal Cohen-Macaulay domain with rational
singularities in characteristic $0$, and $F$-rational in
characteristic $p$.
\end{Corollary}

\begin{proof}
Let $\delta=[a_1\dots a_t\sep b_1\dots b_t]$ and $\gamma=[c_1\dots
c_u\sep d_1\dots d_u]$. Then $\gamma\not\ge\delta$ if $u>t$ or
$c_i<a_i$ or $d_i<b_i$ for some $i=1,\dots,u$. Thus
$\ini(I(X;\delta))$ is generated by those monomials for which
certain exponents are positive. This shows that
$J=\ini(I(X;\delta))$ is a prime ideal.

Therefore the residue class ring $D/J$ is (isomorphic to) a normal
semigroup ring: $D/J$ is a retract of $D$. Now the deformation
arguments apply again.
\end{proof}

Let $\GL=\GL(r,K)$ be the general linear group of invertible $r\times r$-matrices
with entries in $K$. For $f(Y,Z) \in K[Y,Z]$ and $T \in \GL$ we set
$T(f)=f(YT^{-1},TZ)$. This defines a group action on $K[Y,Z]$
as a group of $K$-automorphisms on $K[Y,Z]$.
It turns out that if $|K|= \infty$, then $K[YZ]\cong R_{r+1}$ is the
ring of invariants $K[Y,Z]^{\GL}$ under the action of $\GL$. In the general case
one can show that $K[YZ]$ is the ring of the so-called absolute $\GL$-invariants.

Similar one can consider the action of the special linear group
$\SL=\SL(r,K)=\{T \in \GL(r,K)\colon \det(T)=1\}$ on $K[X,Y]$.
In this case the ring of (absolute) $\SL$-invariants is the $K$-subalgebra
$\Tilde{R}_{r+1}\subset K[Y,Z]$ generated by the entries of $YZ$,
the $r$-minors of $Y$ and the $r$-minors of $Z$.
(See \cite{BV}, Section 7 for definitions and proofs.)
We can study the ring $\Tilde{R}_{r+1}$ analogously to $R_{r+1}$.

\begin{Theorem}\label{invariant}\leavevmode\nopagebreak
\begin{itemize}
\item[(a)] The initial algebra
$\Tilde{D}_{r+1}=\ini(\Tilde{R}_{r+1})\subset K[Y,Z]$
is generated by the monomials
\begin{enumerate}
\item
$Y_{a_11}\cdots Y_{a_tt} Z_{1b_1}\cdots Z_{tb_t}$
with $1\leq t < r$, $a_1<\dots< a_t$ and $b_1<\cdots<b_t$,
\item
$Y_{a_11}\cdots Y_{a_rr}$
with $a_1<\dots< a_r$,
\item
$Z_{1b_1}\cdots Z_{rb_r}$
with $b_1<\dots< b_r$.
\end{enumerate}
\item[(b)] $\Tilde{D}_{r+1}$ is a normal semigroup ring.
\item[(c)] $\Tilde{R}_{r+1}$
is a normal domain,
Cohen-Macaulay, with
rational singularities in characteristic $0$, and $F$-rational in
characteristic $p>0$.
\end{itemize}
\end{Theorem}
\begin{proof}
Let $\pp$ (resp. $\qq$) be the ideal of $K[YZ]\cong R_{r+1}$
generated by the set $\Gamma_r$ (resp.\ $\Gamma_c$) consisting of all
$r$-minors of the first $r$ rows (resp.\ the
first $r$ columns) of the matrix $YZ$. We investigate the ideals
$\pp^t$ and $\qq^t$. The set of all standard bitableaux, which
contain at least $t$ factors of $\Gamma_r$ (resp.\ $\Gamma_c$)
form a $K$-basis of $\pp^t$ (resp. $\qq^t$). (This follows
directly from the fact that $\pp$ and $\qq$ are
straightening-closed ideals of $K[YZ]$; compare \cite[9.6]{BV}.)

$K[Y,Z]$ is a bigraded $K$-algebra in which all entries
of $Y$ have bidegree $(1,0)$ and all entries of $Z$ have bidegree
$(0,1)$. Note that $\Tilde{R}_{r+1}$ is a graded $K$-subalgebra of
$K[Y,Z]$ where $(\Tilde{R}_{r+1})_t$ contains the bihomogeneous
elements $(d_1,d_2)$ such that $d_2-d_1=tr$. In \cite[9.21]{BV} it
is shown that $(\Tilde{R}_{r+1})_t$ is isomorphic to $\pp^t$ as a
$K$-vector space if $t\geq 0$ and isomorphic to $\qq^{-t}$ as a
$K$-vector space if $t\leq0$. This isomorphism is induced by
$$
[a_1\ldots a_r]_Y \mapsto [a_1\ldots a_r\sep 1\ldots r]_{YZ},\quad
[b_1\ldots b_r]_Z \mapsto [1\ldots r\sep b_1\ldots b_r]_{YZ};
$$
Observe that $[a_1\ldots a_r]_Y[b_1\ldots b_r]_Z= [a_1\ldots
a_r\sep b_1 \ldots b_r]_{YZ}$. Then a $K$-basis of
$\Tilde{R}_{r+1}$ consists of the monomials
$$
\prod_{i=1}^{t_1}[a_{i1}\ldots a_{ir}]_Y \cdot \Sigma_1,\quad
\prod_{i=1}^{t_2}[b_{i1}\ldots b_{ir}]_Z \cdot \Sigma_2
$$
where $\Sigma_1,\Sigma_2$ are standard monomials in $K[YZ]\cong R_{r+1}$
and
$$
\prod_{i=1}^{t_1}[a_{i1}\ldots a_{ir}\sep 1\ldots r]_{YZ} \cdot \Sigma_1,\quad
\prod_{i=1}^{t_2}[1 \ldots r \sep b_{i1} \ldots b_{ir}]_{YZ} \cdot \Sigma_2
$$
are standard monomials in $\pp^{t_2}$ (resp.\ $\qq^{t_1}$).
It follows from \ref{iniSt} and the observation before
that the initial monomials are
\begin{align*}
\ini
\Biggl(
\prod_{i=1}^{t_1}[a_{i1},\ldots,a_{ir}]_Y \cdot \Sigma_1
\Biggr) 
&= 
\prod_{i=1}^{t_1}Y_{a_{i1}1}\cdots Y_{a_{ir}r}
\cdot
\ini(\Sigma_1),\\
\ini
\Biggl(
\prod_{i=1}^{t_2}[b_{i1},\ldots,b_{ir}]_Z \cdot \Sigma_2
\Biggr) 
&= \prod_{i=1}^{t_2}Z_{1b_{i1}}\cdots Z_{rb_{ir}}
\cdot \ini(\Sigma_2).
\end{align*}
These distinct monomials are a $K$-basis of $\Tilde{D}_{r+1}$, since
the Hilbert functions of $\Tilde{D}_{r+1}$ and $\Tilde{R}_{r+1}$
coincide. This already proves (a).

To prove (b) one argues similar to the proof of \ref{main}, and (c)
follows again from general theorems on flat deformation.
\end{proof}

\begin{Remark}
\label{equation_2} Again we can describe the set $\Tilde{E}$ of
vectors $[(\alpha_{ij}),(\beta_{uv})]\in
(\RR^{mr})\oplus(\RR^{rn})$ that appear as exponent vectors of
elements in $\Tilde{D}_{r+1}=\ini(\Tilde{R}_{r+1})$. It is the set
of lattice points in the cone defined by the conditions
(\ref{diag1})--(\ref{nonneg}) and
$$
\sum_{i=1}^n \alpha_{ij}-\sum_{v=1}^n\beta_{jv}=0,\qquad
j=1,\dots,r-1.
$$
Note that we have left out exactly one equation from
(\ref{koppel}), namely that for $j=r$.
\end{Remark}

\begin{Remark}
(a)
The program by which Theorems \ref{main} and \ref{invariant} have
been proved consists of three steps: (1) determine the initial
algebra $\ini(R)$ of an algebra $R$ (with respect to a suitable
embedding of $R$ into a polynomial ring), (2) show that $\ini(R)$
is normal and (3) conclude that $R$ is normal, Cohen-Macaulay,
with rational singularities in characteristic $0$, and
$F$-rational in characteristic $p>0$.

This program can also be carried out for several objects derived
from or similar to the rings $R_{r+1}$:
\begin{enumerate}
\item
The Rees algebra $\Dirsum_k \overline I_{s+1}^kT^k \subset
R_{r+1}[T]$ where $s<r$ and $\overline I_{s+1}$ is the ideal
generated by the residue classes of the $s+1$-minors in
$R_{r+1}$.
\item
The symbolic Rees algebra $\Dirsum_k \overline I_{s+1}^{(k)} T^k
\subset R_{r+1}[T]$ where $\overline I_{s+1}^{(k)}$ is the
symbolic powers of $\overline I_{s+1}$.
\item
The subalgebra $A_{r+1,t}$ of $R_{r+1}$ which is generated by the
residue classes of all $t$-minors of the matrix $X$.
\end{enumerate}
For (i) and (iii) one needs that the characteristic of $K$ is $0$
or $>\min\bigl(s+1,m-(s+1),n-(s+1)\bigr)$; see \cite[Section
10]{BV} or \cite{BC4}.

(b) One can also consider a symmetric $n\times n$-matrix $X^\sym$
of indeterminates, i.e.\ $X^\sym_{ij}=X^\sym_{ji}$. In this
situation we have to replace the generic point $K[YZ]$ of $K[X]$
above with the generic point $K[YY^\tr]$ of $K[X^\sym]$ where $Y$
is an $n\times r$ matrix and $Y^\tr$ is the transpose of $Y$. The
proofs are almost the same with minor modifications.

(c) The method presented above provides a comfortable approach to
the structural properties of the determinantal rings. Despite of
the fact that we use term orders it is not a substitute for the
computation of \emph{Gröbner bases} of the determinantal ideals
\emph{within} $K[X]$, or, more precisely, with respect to the
monoid of monomials of $K[X]$. For this task one has to use other
methods, for example the Knuth-Robinson-Schensted correspondence.
(See \cite{BC4} for details.)
\end{Remark}

%
%
%
\section{Cohen-Macaulay and Ulrich ideals}

Suppose that $1\leq r<\min\{m,n\}$. Let $\pp$ and $\qq$ be the
ideals in $R_{r+1}$ as defined in the proof of Theorem
\ref{invariant}: $\pp$ is generated by the $r$-minors of the first
$r$ rows and $\qq$ is generated by the $r$-minors of the first $r$
columns.

Let $J$ be a reflexive rank $1$ module. Then $J$ is isomorphic to
a divisorial ideal.
It is known that the classes $[\pp], [\qq]\in\Cl(R_{r+1})$ are
inverse to each other and that each of them generates (the
infinite cyclic group) $\Cl(R_{r+1})$; e.g.\ see \cite[(8.4)]{BV}.
This implies that all divisorial ideals are represented by the
symbolic powers $\pp^{(t)}$ and $\qq^{(t)}$, $t\ge 0$. Moreover,
$\pp^{(t)}=\pp^t$ and $\qq^{(t)}=\qq^t$ for all $t$
\cite[(9.18)]{BV}.
Thus $J\cong \pp^t$ or $J\cong \qq^t$ for some $t\ge 0$.
Hence, up to isomorphism the powers $\pp^t$ and $\qq^t$
represent all reflexive rank $1$ modules. In this section we study
their Cohen-Macaulay and Ulrich property.

We briefly recall the definition of an Ulrich ideal: Let $S$ be a
homogeneous Cohen-Macaulay $K$-algebra and let $M$ be a finitely
generated graded maximal Cohen-Macaulay $S$-module. Then
$\mu(M)\leq e(M)$ where $\mu(M)$ denotes the minimal number of
generators of $M$ and $e(M)$ denotes the multiplicity of $M$
(e.g.\ see Brennan, Herzog and Ulrich \cite{BHU}). In case of
equality, $M$ is called an \emph{Ulrich module}. A graded ideal
$I\subset S$ is said to be an \emph{Ulrich ideal} if it is an
Ulrich module. If $S$ is a domain and $I\neq 0$ then $e(I)=e(S)$
and hence $I$ is an Ulrich ideal if and only if it is
Cohen-Macaulay and $\mu(I)=e(S)$.\medskip

We start by computing the minimal number of generators for the
powers of the ideals $\pp$ and $\qq$.

\begin{Proposition}
For any integer $t\geq 1$ the number $\mu(\pp^t)$ is equal to the
determi\-nant of the matrix
$$
\left[{t+n-j}\choose{n-i}\right]_{1\leq i,j\leq r}
$$
and the number $\mu(\qq^t)$ is equal to the determinant of the
matrix
$$
\hskip 11pt\left[{t+m-j}\choose{m-i}\right]_{1\leq i,j\leq r}\ .
$$
\end{Proposition}
\begin{proof}
By symmetry it is enough to prove the assertion for $\pp^t$.
According to \cite[(9.3)]{BV} the ideal $\pp^t$ is generated by
the standard bitableaux which are products of exactly $t$
$r$-minors of the first $r$ rows of $X$ (modulo $I_{r+1}$). These
standard bitableaux are $K$-linearly independent.
Their number coincides with the number of standard bitableaux with $t$ factors
in the coordinate ring $G(r,n)$ of the Grassmannian of
$r$-dimensional vector spaces in $K^n$ because
the latter elements are the preimages  of the generators of
$\pp^t$  in $K[X]$. So we can finish our proof
quoting the classical formula of Hodge (for example, see Ghorpade
\cite[Theorem 6]{G}) by which $\dim_K
G(r,n)_{t}$ is equal to the determinant of the matrix given in the
assertion.
\end{proof}
Next we show:

\begin{Lemma}
The multiplicity of $R_{r+1}$ coincides with $\mu(\pp^{m-r})$ and $\mu(\qq^{n-r})$.
\end{Lemma}
\begin{proof}
The multiplicity of $R_{r+1}$ is known to be the determinant of the matrix
$$
B=\left[{m+n-i-j}\choose{n-j}\right]_{1\leq i,j\leq r}.
$$
(E.g.\ see Herzog and Trung \cite{HT}.) By the above proposition
we know that $\mu(\pp^{m-r})$ is equal to the determinant of the
matrix
$$
\hskip 10pt A=\left[{m+n-r-j}\choose{n-i}\right]_{1\leq i,j\leq r}\ .
$$
Using the binomial identity
${{a}\choose{b}}+{{a}\choose{b+1}}={{a+1}\choose{b+1}}$ one can
transform $A$ into the transpose of $B$ by elementary row
operations\rule{0pt}{11pt} which do not affect the determinant.
This proves $\mu(\pp^{m-r})=e(R_{r+1})$. The equation
$\mu(\qq^{n-r})=e(R_{r+1})$ can be obtained in an analogous way.
\end{proof}

As a function of $t$ the minimal number of generators
$\mu(\pp^{t})$ is evidently a strictly increasing function in $t$.
Thus $\mu(\pp^{t})>\mu(\pp^{m-r})=e(R_{r+1})$ for $t>m-r$, and
$\pp^{t}$ cannot be a Cohen-Macaulay ideal. By the same reason
$\qq^{t}$ cannot be Cohen-Macaulay for $t>n-r$.

\begin{Theorem}\label{CMUlrich}
Let $\hskip 1pt t\geq 1$ be an integer.  The power $\pp^t$ (resp.\
$\qq^t$) is a Cohen-Macaulay ideal if and only if $t\leq m-r$
(resp.\ $t\leq n-r$). The powers $\pp^{m-r}$ and $\qq^{n-r}$ are
both Ulrich ideals.
\end{Theorem}

\begin{proof}
The crucial point which has not been proved yet is that $\pp^t$
(resp.\ $\qq^t$) is a Cohen-Macaulay ideal for $t\leq m-r$ (resp.\
$t\leq n-r$). By symmetry it is enough to deal with $\pp^t$.

Assume that $t\leq m-r$.
We consider the set of all standard bitableaux of $R_{r+1}$,
which contain at least $t$ factors of the generators of $\pp$.
We already observed in the proof of
\ref{invariant} that these ele\-ments form
a $K$-basis of $\pp^t$.
Now we use the generic point $\phi:R_{r+1} \to K[Y,Z]$ to embed
$\pp^t$ into $K[Y,Z]$, and investigate the initial ideal
$\aaf_t=\ini(\phi(\pp^t))\subset D_{r+1}$.

Let $E_t$ be the subset of $E$ (compare \ref{equation} for the
definition of $E$) consisting of all vectors in
$(\RR^{mr})\oplus(\RR^{rn})$ that appear as
exponent vectors of the elements in $\aaf_t$. One easily checks
that
\begin{align*}
E_t& =\{\,[(\alpha_{ij}),(\beta_{uv})]\in E\mid \alpha_{ii}\geq t, \enskip i=1,\ldots,r \,\} \\
   & =\{\,[(\alpha_{ij}),(\beta_{uv})]\in E\mid \alpha_{rr}\geq t\,\}.
\end{align*}
We want to show that $\aaf_t$ is a conic ideal in $D_{r+1}$ (see
Bruns and Gubeladze \cite[Section 3]{BG1}). To this end we have to
find $w_t\in\RR E$ such that $E_t=\ZZ E\cap(w_t+\RR_{+}E)$. Note
that $\RR E$ is the set of all vectors
$[(\alpha_{ij}),(\beta_{uv})]\in(\RR^{mr})\oplus(\RR^{rn})$
that satisfy the equations
$$
\alpha_{ij}=\beta_{uv}=0,\quad j>i, u>v, \qquad
\sum_{i=1}^n\alpha_{ij}-\sum_{v=1}^n\beta_{jv}=0,\quad
j=1,\ldots,r\,,
$$
and that $\ZZ E=\RR E \cap (( \ZZ^{mr})\oplus(\ZZ^{rn} ))$. We choose a
positive real number $\varepsilon<1$ and define
$\rule{0pt}{12pt}w_t=[(\alpha_{ij}),(\beta_{uv})]$ by setting
$$
\alpha_{ij}=\begin{cases}
t-\epsilon, & \text{if }\ i=j,\\
\rule{0pt}{12pt}-(t-\epsilon)/(m-r), & \text{if }\ j<i\leq m-r+j,\\
\rule{0pt}{12pt}0, & \text{otherwise}.
\end{cases}
$$
and $\beta_{uv}=0$ for all $u,v$. It is clear that $w_t\in\RR E$.
Since $-(t-\varepsilon)/(m-r)>-1$ (this is the point where we need
$t\leq m-r\,$!) we have $\ZZ E\cap(w_t+\RR_{+}E)=E_t$. So $\aaf_t$
is indeed a conic ideal. Since every conic ideal in a normal
semigroup ring is Cohen-Macaulay (see \cite[3.3]{BG1}) we
conclude that $\aaf_t$ is a Cohen-Macaulay ideal in the ring
$D_{r+1}$. But this implies that $\pp^t$ is a Cohen-Macaulay ideal
in the ring $R_{r+1}$ (e.g.\ see \cite[3.16]{BC4}).
\end{proof}

The case $r=1$ of the theorem has been proved (and the general has
been conjectured) by Bruns and Guerrieri \cite{BrGu}.

\begin{Corollary}
The ideals $\pp^t$, $0\leq t\leq m-r$, and $\qq^t$, $0<t\leq n-r,$
represent all isomorphism classes of maximal Cohen-Macaulay
$R_{r+1}$-modules of rank $1$.
\end{Corollary}

\begin{proof}
Let $M$ be a maximal Cohen-Macaulay $R_{r+1}$-module of rank 1.
Then $M$ is torsionfree and therefore it is isomorphic to a
fractionary ideal $J$ of $R_{r+1}$. Using the reflexivity
criterion of \cite[1.4.1]{BH}, one sees that $J$ is reflexive and
hence divisorial.

We already noticed in the beginning of this section
that then $J\cong \pp^t$ or $J\cong \qq^t$ for some $t\ge 0$, and the
corollary follows immediately from Theorem \ref{CMUlrich}.
\end{proof}

%
%
%

\end{document}